\documentclass[11pt,reqno]{amsart}

\usepackage{amssymb, amsmath}
\usepackage{hyperref}
\usepackage[nohug]{diagrams}
\diagramstyle[labelstyle=\scriptstyle]

\newtheorem{proposition}{Proposition}[section]
\newtheorem{lemma}[proposition]{Lemma}

\newtheorem{theorem}[proposition]{Theorem}

\newtheorem{problem}{Problem}[]

\theoremstyle{definition}
\newtheorem{definition}[proposition]{Definition}
\newtheorem{example}[proposition]{Example}

\newtheorem{remark}[proposition]{Remark}

\newcommand{\thlabel}[1]{\label{th:#1}}
\newcommand{\thref}[1]{Theorem~\ref{th:#1}}
\newcommand{\selabel}[1]{\label{se:#1}}
\newcommand{\seref}[1]{Section~\ref{se:#1}}
\newcommand{\lelabel}[1]{\label{le:#1}}
\newcommand{\leref}[1]{Lemma~\ref{le:#1}}
\newcommand{\prlabel}[1]{\label{pr:#1}}
\newcommand{\prref}[1]{Proposition~\ref{pr:#1}}

\newcommand{\relabel}[1]{\label{re:#1}}
\newcommand{\reref}[1]{Remark~\ref{re:#1}}
\newcommand{\exlabel}[1]{\label{ex:#1}}

\newcommand{\delabel}[1]{\label{de:#1}}
\newcommand{\deref}[1]{Definition~\ref{de:#1}}

\newcommand{\analog}{\rule{1cm}{0.01mm} }

\begin{document}

\title{When is the diagonal functor Frobenius?}
\author{Alexandru Chirv\u asitu}

\address{Faculty of Mathematics and Computer Science, University of Bucharest, Str.
Academiei 14, RO-70109 Bucharest 1, Romania}
\email{chirvasitua@gmail.com}


\begin{abstract}

Given a complete, cocomplete category $\mathcal C$, we investigate the problem of describing those small categories $I$ such that the diagonal functor $\Delta:\mathcal C\to {\rm Functors}(I,\mathcal C)$ is a Frobenius functor. This condition can be rephrased by saying that the limits and the colimits of functors $I\to\mathcal C$ are naturally isomorphic. We find necessary conditions on $I$ for a certain class of categories $\mathcal C$, and, as an application, we give both necessary and sufficient conditions in the two special cases $\mathcal C={\bf Set}$ or $_R\mathcal M$, the category of left modules over a ring $R$. 

\end{abstract}

\maketitle

\section*{Introduction}

Functors having a left adjoint which is also a right adjoint were investigated by Morita in \cite{Mo}, where it is shown that given a ring morphism $R\to S$, the restriction of scalars functor has this property if and only if $R\to S$ is a Frobenius extension: $S$ is finitely generated and projective in $_R\mathcal M$, and $S\cong\ _R{\rm Hom}(S,R)$ as $(S,R)$-bimodules. Pairs of functors $F,G$ (between module categories) with the property that both $(F,G)$ and $(G,F)$ are adjunctions are called by Morita {\it strongly adjoint pairs of functors}. Later, a functor $F$ having a left adjoint which is also a right adjoint came to be referred to as a {\it Frobenius functor} (\cite{CMZ}), and Morita's strongly adjoint pairs of functors are now known as {\it Frobenius pairs}. 

The natural question arises of when various well-known and extensively used functors are Frobenius. Examples include the already mentioned case of the restrictions of scalars functor for a ring extension (\cite{Men, Mo}), forgetful functor from Doi-Hopf (or Doi-Koppinen) modules to modules (\cite{CMZ}), forgetful functor from $G$-graded modules over a $G$-graded ring to modules, where $G$ is a group (\cite{DasNas}), corestriction of scalars through an $A$-coring map $C\to D$ (\cite{Iov2}, or \cite{ZD} in the more general setting of a map from an $A$-coring $C$ to a $B$-coring $D$), and many more. 

In this paper the point of view is the following one: we fix a complete, cocomplete category $\mathcal C$, and seek to characterize those small categories $I$ for which the functors $\mathcal C^I\to C$ sending a functor in $\mathcal C^I$ to its limit and colimit are naturally isomorphic. We call such a category $\mathcal C$-Frobenius. The connection to Frobenius functors (hence the name $\mathcal C$-Frobenius) is highlighted by the following observation: the functor $\varprojlim:\mathcal C^I\to\mathcal C$ is right adjoint to the diagonal functor $\Delta:\mathcal C\to\mathcal C^I$, whereas the colimit functor is the left adjoint to $\Delta$. Hence our question can be rephrased as follows: for which small categories $I$ (depending on $\mathcal C$) is the diagonal functor $\Delta:\mathcal C\to\mathcal C^I$ a Frobenius functor? 

This question is investigated in \cite{Iov1}, for discrete small categories $I$ (i.e. sets), and categories $\mathcal C$ enriched over the category of commutative monoids (referred to as ${\bf AMon}$ categories), and having a zero object. In that setting the problem is to find those sets $I$ for which direct sums and direct products in $\mathcal C$ indexed by $I$ are naturally isomorphic. The main result \cite[Proposition 1.3]{Iov1} says that under reasonably mild conditions, this is equivalent to $I$ being finite. 

Here, on the other hand, we focus mainly on connected categories $I$. The structure of the paper is as follows:

In \seref{1} we introduce some conventions and prove \leref{discr_conn}, which allows us later on to break up the main problem into the two cases when $I$ is discrete (a set) or connected. 

In \seref{2} we introduce the class of categories $\mathcal C$ we will be concerned with, which we call {\it admissible}, and also turn our attention to the case when $I$ is connected. Two general results, \thref{Frob>inv} and \prref{cat<>grp}, are proven in this setting. 

In \seref{3} necessary and sufficient conditions on $I$ are found in order that it be ${\bf Set}$-Frobenius or $_R\mathcal M$-Frobenius, where $R$ is a ring and $_R\mathcal M$ is the category of left $R$-modules. Since both ${\bf Set}$ and $_R\mathcal M$ are admissible in the sense of \seref{2}, the results proven there can be applied to the two particular cases. 

The conditions on $I$ appearing in the main results of \seref{3} (Theorems 3.1 and 3.2) are of a combinatorial nature. The full description of the statements of these theorems requires some preparation (\deref{inv}), but they immediately imply, for instance, the characterization of ${\bf Set}$ or $_R\mathcal M$-Frobenius monoids $I$ (as usual, we regard a monoid as a category with a single object). A consequence of \thref{Set} is that the ${\bf Set}$-Frobenius monoids $I$ are precisely those containing an element $a\in I$ which is a fixed point for all left and right multiplications: $xa=ax=a,\ \forall x\in I$. Similarly, \thref{R-Mod} implies that a monoid $I$ is $_R\mathcal M$-Frobenius if and only if it contains a finite (non-empty) set $S$ on which all multiplications, left or right, act as permutations, and such that the cardinality $|S|$ of $S$ is invertible in the ring $R$. The full description of connected Frobenius categories $I$ in the two cases is a natural generalization of this discussion. 

Finally, in \seref{4} we finish with some open problems for the reader. 

\section{Preliminaries}\selabel{1}

Throughout this paper, $\mathcal C$ will denote a complete, cocomplete category, while $I$ stands for a small category. In general, for notions pertaining to category theory, we refer to \cite{MacL}. The convention for composing morphisms is the usual one: given two morphisms $f:x\to y$ and $g:y\to z$ in a category, their composition is $gf:x\to z$. In order to keep the notation simple, if $i$ is an object of $I$ we write $i\in I$ (rather than $i\in{\rm Ob}(I)$, for example). Sometimes, in order to make it easier to keep track of the objects involved in morphisms, we shall denote $f\in {\rm Hom}(i,j)$ by $f_i^j$. Similarly, we might denote a subset $S\subseteq {\rm Hom}(i,j)$ by $S_i^j$. Given a set $S\subseteq {\rm Hom}(i,j)$ and a morphism $f\in {\rm Hom}(j,k)$, $fS$ stands for the set of all morphisms $fg,\ g\in S$; similarly for $Sf$, when the composition makes sense. Given categories $X,Y$, we denote the category ${\rm Functors}(X,Y)$ simply by $Y^X$. All functors are covariant, except when explicitly mentioned otherwise.

\begin{definition}

Let $\mathcal C$ be a complete, cocomplete category. A small category $I$ is said to be {\it $\mathcal C$-Frobenius} if the diagonal functor $\Delta:\mathcal C\to\mathcal C^I$ is a Frobenius functor.

\end{definition}

\begin{remark}

The left adjoint to $\Delta$ is the functor $\mathcal C^I\to\mathcal C$, sending $F\in\mathcal C^I$ to its colimit $\varinjlim F$. Similarly, the right adjoint to $\Delta$ is the functor sending $F\in\mathcal C^I$ to its limit $\varprojlim F$ (\cite[Chapter IV $\S$2]{MacL}). Consequently, saying that $\Delta$ is Frobenius is the same as saying that $\varprojlim$ and $\varinjlim$ are naturally isomorphic. This means that we can find, for each functor $F\in\mathcal C^I$, an isomorphism $\psi_F:\varprojlim F\to\varinjlim F$ such that for every natural transformation $\eta:F\to G$ one has the commutative diagram

$$
\begin{diagram}[labelstyle=\scriptstyle]
\varprojlim F             &\rTo^{\psi_F}              &\varinjlim F\\
\dTo^{\varprojlim\eta}    &                           &\dTo_{\varinjlim\eta}\\
\varprojlim G             &\rTo^{\psi_G}              &\varinjlim G
\end{diagram}
$$

\end{remark}

\begin{remark}

Notice that the empty category is $\mathcal C$-Frobenius if and only if $\mathcal C$ has a zero object. In order to avoid splitting the arguments into cases, we assume from now on that all our categories are non-empty. 

\end{remark}

We remarked earlier that we would be concerned primarily with the case when $I$ is connected. In fact, as the following lemma shows, the general problem of finding the $\mathcal C$-Frobenius small categories $I$ for a given $\mathcal C$ breaks up into the connected and the discrete case under certain conditions which do occur in the cases of interest.

\begin{lemma}\lelabel{discr_conn}

Let $\mathcal C$ be a complete, cocomplete category and $I$ a small category with connected components $I_j,\ j\in J$. Then:

\begin{enumerate}
\renewcommand{\labelenumi}{(\alph{enumi})}

\item If each component $I_j$ is $\mathcal C$-Frobenius and the set $J$, regarded as a discrete category, is $\mathcal C$-Frobenius, then $I$ is $\mathcal C$-Frobenius.

\item If $I$ is $\mathcal C$-Frobenius, then $J$ is $\mathcal C$-Frobenius.

\item The converse of ${\rm (a)}$ holds if $\mathcal C$ has a zero object.

\end{enumerate}

\end{lemma}

\begin{proof}

Before proving the three assertions, we make some observations useful in all three arguments. Fix a functor $F\in\mathcal C^I$, and consider the contravariant functor $T_F:\mathcal C\to{\bf Set}$ defined by sending each object $c$ to the set of cones $\tau:c\stackrel{\cdot}{\rightarrow}F$ (Mac Lane's terminology and notation; see \cite[Chapter III $\S$3]{MacL}). The set of cones can also be defined as the object set of the comma category $c\downarrow F$ (\cite[Chapter II $\S$6]{MacL}). Since $I$ is small, the comma category is indeed small, so it makes sense to talk about its object set. Notice that the limit $\varprojlim F$ is precisely the representing object of $T_F$. Moreover, $F\mapsto T_F$ is natural in $F$. 

On the other hand, again having fixed $F\in\mathcal C^I$, consider the functor $S_F:\mathcal C\to{\bf Set}$ sending $c$ to collections of cones $\tau_j:c\stackrel{\cdot}{\rightarrow}(F|_{I_j}),\ j\in J$ from $c$ to the restrictions of $F$ to the connected components $I_j$. By the definition of limits, the representing object for $S_F$ is $\displaystyle\prod_{j\in J}(\varprojlim F|_{I_j})$. Notice however that, since there are no morphisms between distinct connected components, the functors $T_F$ and $S_F$ actually coincide. In conclusion, the representing objects $\varprojlim F$ and $\displaystyle\prod_{j\in J}(\varprojlim F|_{I_j})$ are in fact isomorphic; the isomorphism exhibited here is natural in $F$, because $F\mapsto T_F$ is. Similarly, $\varinjlim F\cong\displaystyle\coprod_{j\in J}(\varinjlim F|_{I_j})$. We are now ready for the proof proper.

(a) We have just seen that $\varprojlim F\cong\displaystyle\prod_{j\in J}(\varprojlim F|_{I_j})$ naturally in $F$. Each $I_j$ is Frobenius, so the latter is isomorphic to $\displaystyle\prod_{j\in J}(\varinjlim F|_{I_j})$ (naturally in $F$). The component set $J$ is Frobenius, so this is isomorphic to $\displaystyle\coprod_{j\in J}(\varinjlim F|_{I_j})$ (again, naturally in $F$). Finally, the above discussion shows that this is isomorphic to $\varinjlim F$.

(b) Instead of looking at the whole of $\mathcal C^I$, consider only those functors $I\to\mathcal C$ which restrict to constants on each component $I_j$. These are precisely the functors factoring through the canonical functor $\nu:I\to J$, which sends each $I_j$ to $j$. Again, use the isomorphism $\varprojlim F\cong\displaystyle\prod(\varprojlim F|_{I_j})$: the limit of a constant functor on a connected category is easily seen to be precisely the image object (with all structural morphisms equal to the identity); it follows that in the case at hand, when $F$ restricts to a constant on each component, $\varprojlim F$ is naturally isomorphic to the product of the objects $F(I_j)$. The same discussion applies to colimits: $\displaystyle\varinjlim F\cong\coprod F(I_j)$. The desired conclusion that $J$ must be $\mathcal C$-Frobenius follows.

(c) In view of (b), we must show that given the additional hypothesis of a zero object, each $I_j$ is $\mathcal C$-Frobenius. Fix some index $k\in J$, and consider only those functors $I\to\mathcal C$ which send each component $I_j,\ j\ne k$ to the zero object $0$. Using once more the discussion at the beginning of the proof, we conclude that for these functors, the limit is naturally isomorphic to the product $\displaystyle\left(\varprojlim F|_{I_k}\right)\times\prod_{j\ne k}0$. Since in any complete category product with the final object is naturally isomorphic to the identity, we conclude that $\varprojlim F\cong\varprojlim F|_{I_k}$, naturally in $F$. Similarly, the colimit of $F$ is isomorphic to that of $F|_{I_k}$, so $I_k$ must indeed be $\mathcal C$-Frobenius if $I$ is.  
\end{proof}

\section{Admissible categories, free objects, and some general results}\selabel{2}

In the end, we are going to find the small categories $I$ which are ${\bf Set}$-Frobenius and those which are $_R\mathcal M$-Frobenius for a given ring $R$. Part of that proof will be unified by the results in this section, dealing with a certain class of categories $\mathcal C$ which contains both ${\bf Set}$ and $_R\mathcal M$, and many more familiar categories. We introduce this class below:

\begin{definition}\delabel{adm}

A category $\mathcal C$ is called {\it admissible} if:

\begin{enumerate}

\item it is complete and cocomplete
\item there is a faithful functor $U: \mathcal C\to{\bf Set}$ which has a left adjoint $T$
\item for at least one object $c$ of $\mathcal C$, the set $Uc$ has $\ge 2$ elements
\item for any set $X$ and any element $t$ of the set $UT(X)$, there is a smallest finite subset $Y\subseteq X$ such that $t$ belongs to the set $UT(Y)$.

\end{enumerate}

\end{definition}

For a set $X$, we denote the free object $T(X)$ by $T_X$. The faithful functor $U:\mathcal C\to {\bf Set}$ makes $\mathcal C$ into what is usually called a {\it concrete category}. Most of the time we will simply omit $U$, and regard $\mathcal C$ as a category whose objects are sets (with ``additional structure''; that is, we keep in mind that the same set might correspond to different objects), and whose morphisms are functions between these sets.

\begin{remark}\relabel{defrem1}

Condition (3) implies that for each set $X$, the component $\psi_X:X\to UT_X$ of the unit of our adjunction $(T,U)$ is mono. Indeed, if $c$ is an object of $\mathcal C$ such that the set $Uc$ has at least two elements and $X$ is any set, then any two different elements of $X$ can be mapped to different elements of $Uc$, meaning that any two different elements of $X$ must have different images in the set $UT_X$. Hence, from now on we will regard $X$ as a subset of $UT_X$ (or of $T_X$, with the convention in the previous paragraph). Also, condition (3) implies that $T_\emptyset$ is not isomorphic to any of the other free objects, a fact that will be useful at some point: $T_\emptyset$ is initial, whereas any other free object admits at least two morphisms to any object $c\in\mathcal C$ such that $Uc$ has at least two elements.

\end{remark}

\begin{remark}\relabel{defrem2}

Another observation which will be used tacitly from now on is this: inclusions of sets $X\to Y$ induce inclusions of sets $T_X\to T_Y$ (we omit $U$ in this remark). When $X\ne\emptyset$ this is clear: every monomorphism of sets $X\to Y$ is then a coretraction, and functors preserve coretractions. When $X=\emptyset$, on the other hand, $T_X$ is the initial object of $\mathcal C$. The initial object can be constructed, in any complete category, as a subobject of any weakly initial object (see \cite[Chapter V $\S$6, proof of Theorem 1]{MacL}). More precisely, it is the equalizer of all endomorphisms of any such object. By weakly initial we mean object admitting a morphism (not necessarily unique) to any object. Free objects are all weakly initial (unless $T_\emptyset=\emptyset$, in which case there is nothing left to prove), so $T_\emptyset$ is a subobject of each of them. Right adjoints (such as $U$) preserve monomorphisms, so, given a subset $X$ of $Y$, we will regard $T_X$ as a subset of $T_Y$; the inclusion is always the one induced by $X\to Y$. 

\end{remark}

Here we make a short digression to identify many familiar categories which are, in fact, admissible. These are the so-called {\it varieties of algebras}, in the sense of Universal Algebra. For definitions and a detailed treatment we refer to \cite[Chapter II]{BurrSank}. Also, there is some discussion on the topic, from a more category theoretical point of view, in \cite[Chapter V $\S$6]{MacL}; here the main definitions are given, and the proof for the existence of free objects is sketched, using Freyd's Adjoint Functor Theorem (\cite[Chapter V $\S$6, Theorem 2]{MacL}). 

We will not give complete proofs or definitions here. Given an $\mathbb N$-graded set $\Omega$ whose elements are called {\it operations}, an action of $\Omega$ on a set $A$ is a map assigning to each $\omega\in\Omega$ of degree $n\in\mathbb N$ a function $\omega_A:A^n\to A$. The degree $n$ is also called the {\it arity} of $\omega$. From the operations in $\Omega$, named fundamental operations, others can be derived, by composition and substitution; see the reference from Mac Lane. A set $E$ of {\it equational identities} is a set of pairs $(\mu,\nu)$ of derived operations having the same arity. A set $A$ with an $\Omega$ action is then said to satisfy the equations $E$ if $\mu_A=\nu_A$ for all $(\mu,\nu)\in E$. The class of all sets with an $\Omega$ action and satisfying the identities $E$ will be denoted by $\langle\Omega,E\rangle-{\bf Alg}$, and a member of this class will be called an $\langle\Omega,E\rangle$-algebra. 

A morphism between algebras $A,B\in\langle\Omega,E\rangle-{\bf Alg}$ is a map $f:A\to B$ which, for each $\omega\in\Omega$, makes the following diagram commutative:

$$
\begin{diagram}[labelstyle=\scriptstyle]
A^n                &\rTo^{f^n}               &B^n\\
\dTo^{\omega_A}    &                         &\dTo_{\omega_B}\\
A                  &\rTo^f                   &B
\end{diagram}
$$
We now have a category $\langle\Omega,E\rangle-{\bf Alg}$. Examples include the categories of sets (no operations at all), monoids, groups, rings, modules over a ring $R$ (these are abelian groups with some unary operations describing multiplications with scalars in $R$), Lie algebras, etc. Notice that we allow the underlying set of an algebra to be empty, although the authors of \cite{BurrSank} do not. A variety contains the empty set if and only if there are no nullary operations (i.e. operations of arity $0$). 

The definitions allow for a variety of algebras not to satisfy condition (3) of \deref{adm}. Assuming it does, however, it can be shown that $\langle\Omega,E\rangle-{\bf Alg}$ is admissible, with $U$ (from \deref{adm}) being the forgetful functor, which sends an algebra to its underlying set, and a morphism of algebras to the corresponding map of sets. We will not give the complete proof here. As mentioned above, Mac Lane proves the existence of free objects indirectly, using the Adjoint Functor Theorem. In \cite[Chapter II $\S$10]{BurrSank} an explicit construction of free objects is given. Condition (4) follows from the fact that a filtered union of $\langle\Omega,E\rangle$-algebras is again such with an obvious structure; it is easy to check the required universality property for the union of all $T_Y$ as $Y$ ranges through the finite subsets of $X$, which makes it into the free object on $X$; hence, elements of free objects are contained in finitely generated free subobjects. Finally, (4) is proven by noticing that for varieties of algebras, one always has $T_Y\cap T_Z=T_{Y\cap Z}$ (including the case when there are no nullary operations, and the second set in this equality happens to be empty). This can be seen by constructing the free objects explicitly. For completeness and cocompleteness one can mimic the usual constructions of products, coproducts, equalizers and coequalizers from group theory, for example. 

In particular, ${\bf Set}$ and $_R\mathcal M$ are admissible. Of course, this can be seen directly.

We will need the next lemma in the proof of \thref{Frob>inv}.

Recall that a directed graph (digraph) is said to be {\it strongly connected} if for any two vertices $i,j$ there is a directed path from $i$ to $j$. A digraph is said to be {\it transitive} if whenever we have directed paths $i\to j$ and $j\to k$ we also have a directed path $i\to k$. The underlying graph of a category is transitive, for instance. If a digraph is transitive then strong connectedness is equivalent to having an edge $i\to j$ for any pair of distinct vertices $i,j$.

\begin{lemma}\lelabel{strcon}

Let $\mathcal C$ be an admissible category, and $I$ a small, connected, $\mathcal C$-Frobenius category. Then $I$ is in fact strongly connected, i.e. ${\rm Hom}(i,j)\ne\emptyset$ for all pairs of objects $i,j\in I$. 

\end{lemma}

\begin{proof}

We will make use of the following well-known combinatorial result: if a connected directed graph is not strongly connected, then its vertex set can be partitioned into two non-empty subsets $A,\ B$ such that all the arrows connecting them go from $A$ to $B$. Moreover, $A$ can be chosen to be connected. Assuming that $I$ is not strongly connected, apply this to the underlying graph of $I$. We get non-empty, full subcategories $A,\ B$ of $I$ with $A$ connected, which partition its object set, and such that all morphisms between $A$ and $B$ go from $A$ to $B$.  

Now consider the functor $F\in\mathcal C^I$ which restricts to the constant functor $T_\emptyset$ on $A$, to the constant $T_1$ on $B$, and sends all morphisms $A\to B$ onto the unique morphism $T_\emptyset\to T_1$:

$$
\begin{diagram}[labelstyle=\scriptstyle]
A                &                                        &T_\emptyset\\
\dTo\dTo         &\stackrel{F}{\longrightarrow}           &\dTo\\
B                &                                        &T_1
\end{diagram}
$$

An argument very similar to the one used in the proof of \leref{discr_conn} (the beginning of that proof) shows that the limit of $F$ is $T_\emptyset$. On the other hand, the colimit is the coproduct of one copy of $T_1$ for each connected component $B_j,\ j\in J$ of $B$; here $J$ is simply the (non-empty) set of connected components. $T$ is a left adjoint by definition, so it preserves coproducts; this means that $\displaystyle\coprod_J T_1\cong T_J$. We have already remarked, in the discussion after \deref{adm}, that $T_\emptyset$ cannot be isomorphic to a free object $T_J,\ J\ne\emptyset$, so $I$ is not $\mathcal C$-Frobenius. We have reached a contradiction.
\end{proof}

\begin{remark}\relabel{noninit}

Notice that in the above proof, instead of the unique arrow $T_\emptyset\to T_1$ we could just as well have taken the unique arrow from an initial object to a non-initial object. Hence the statement holds for any (complete, cocomplete) category $\mathcal C$ having at least one object which is not initial.  

\end{remark}

The following definition is crucial in subsequent results. $I$ stands for a small category.

\begin{definition}\delabel{inv}

A {\it left invariant system} (LS) of $I$ is a collection of finite, non-empty sets $S_i^j\subseteq {\rm Hom}(i,j)$, one for each pair $i,j\in I$, such that composition to the left with any $f_j^k\in {\rm Hom}(j,k)$ sends $S_i^j$ bijectively onto $S_i^k$ for all $i,j,k\in I$. 

A {\it right invariant system} (RS) of $I$ is a collection of finite, non-empty sets $S_i^j\subseteq {\rm Hom}(i,j)$, one for each pair $i,j\in I$, such that composition to the right with any $f_k^i\in {\rm Hom}(k,i)$ sends $S_i^j$ bijectively onto $S_k^j$ for all $i,j,k\in I$.

An {\it invariant system} (IS) of $I$ is an LS which is also an RS. 

\end{definition}

The main result of this section follows:

\begin{theorem}\thlabel{Frob>inv}

Let $\mathcal C$ be an admissible category, and let $I$ be a small, connected, $\mathcal C$-Frobenius category. Then $I$ has an IS. 

\end{theorem}

\begin{proof}

The functors in $\mathcal C^I$ we will work with are $i^*=T_{{\rm Hom}(i,-)}$ for objects $i\in I$. $T$ being a left adjoint, it preserves colimits. In other words, $\varinjlim i^*\cong T_{\varinjlim {\rm Hom}(i,-)}$. By the description of colimits in ${\bf Set}$ one sees immediately that $\varinjlim {\rm Hom}(i,-)$ is a singleton. In conclusion, $\varinjlim i^*\cong T_1$. By the $\mathcal C$-Frobenius property we can identify $\varprojlim i^*$ with $T_1$ as well. We will denote by $1$ the element generating $T_1$; in the present context it corresponds to the image of any morphism in ${\rm Hom}(i,j)$ through the canonical map ${\rm Hom}(i,j)\to T_{\varinjlim {\rm Hom}(i,-)}\cong\varprojlim i^*$. 

Let $\psi_i^j:T_1\cong\varprojlim i^*\to T_{{\rm Hom}(i,j)}$ be the structure map of the limit, and denote by $x_i^j$ the element $\psi_i^j(1)\in T_{{\rm Hom}(i,j)}$ (keep in mind the convention made after \deref{adm}: we regard the objects of $\mathcal C$ simply as sets, omitting the faithful functor $U:\mathcal C\to{\bf Set}$). By condition (4) of \deref{adm}, there is a smallest finite set $S\subseteq {\rm Hom}(i,j)$ such that $x_i^j\in T_S$. Denote it by $S_i^j$; as the notation suggests, these will be the components of our IS. 

{\it $(S_i^j)$ is an LS}. For all $j,k\in I$ and all $f_j^k$ we have a commutative diagram

$$
\begin{diagram}[labelstyle=\scriptstyle]
T_1               &                  &\\
\dTo^{\psi_i^j}   &\rdTo^{\psi_i^k}  &\\
i^*(j)            &\rTo^{i^*f_j^k}     &i^*(k) 
\end{diagram}
$$ 
It follows that $(i^*f_j^k)(x_i^j)=x_i^k$, so, by the definition of the sets $S_i^j$, we have $f_j^kS_i^j\supseteq S_i^k$. In other words, composition to the left maps $S_i^j$ onto a set containing $S_i^k$. A consequence of this is that $|S_i^k|\le|S_i^j|$ whenever the hom set ${\rm Hom}(j,k)$ is non-empty. However, we know from \leref{strcon} that all hom sets are nonempty, so all $S_i^j$ have the same cardinality. Moreover, composition to the left with any morphism must be a bijection.

All we need to do now in order to conclude that $S=(S_i^j)$ is an LS is to show that the sets $S_i^j$ are non-empty. Assume they are. Then $\psi_i^j$ maps $\varprojlim i^*\cong T_1$ into $T_\emptyset\subset T_{{\rm Hom}(i,j)}=i^*(j)$ for all $j$. This means that the limiting cone $\varprojlim i^*\stackrel{\cdot}{\rightarrow}i^*$ factors through $T_\emptyset\stackrel{\cdot}{\rightarrow}i^*$ which, in turn, implies that $T_1\cong\varprojlim i^*\cong T_\emptyset$. This is impossible by condition (3) in \deref{adm} (see \reref{defrem1}).

{\it $(S_i^j)$ is an RS}. This is where the naturality of $\eta:\varprojlim\cong\varinjlim$ comes in. More pecisely, consider any morphism $f=f_i^j\in {\rm Hom}(i,j)$. It induces a natural transformation $f^*$ from $j^*$ to $i^*$. The corresponding transformations $\varprojlim j^*\to\varprojlim i^*$ and $\varinjlim j^*\to\varinjlim i^*$ will again be denoted by $f^*$. For each $k\in I$ we have the following commutative diagram:

$$
\begin{diagram}[labelstyle=\scriptstyle]
\varinjlim j^*         &\rTo^{\cong}        &\varprojlim j^*          &\rTo^{\psi_j^k}              &T_{{\rm Hom}(j,k)}\\
\dTo^{f^*}             &                    &\dTo^{f^*}               &                             &\dTo^{f^*}  \\
\varinjlim i^*         &\rTo^{\cong}        &\varprojlim i^*          &\rTo^{\psi_i^k}              &T_{{\rm Hom}(i,k)}                           
\end{diagram}
$$   
The horizontal arrows of the left square are the components of the natural isomorphism $\eta:\varinjlim\cong\varprojlim$. 

Notice that $1\in T_1\cong\varinjlim j^*$ gets mapped onto $1\in T_1\cong\varinjlim i^*$ (see the description of $1$ in the first paragraph of the proof). Since we have identified $\varprojlim j^*$ to $T_1$ through $\eta$, it follows from this diagram that $f^*(x_j^k)=x_i^k$. By the definiton of the sets $S$, this means that $S_j^kf_i^j\supseteq S_i^k$. Now we continue as in the proof for left invariance, using the fact that all hom sets are non-empty ($I$ is strongly connected).  
\end{proof}

Let $I$ be a small, connected category with an IS $(S_i^j)$ (in particular, $I$ will be strongly connected). Consider a set $S_i^i$ for some object $i\in I$. Composition of morphisms gives such a set a structure of finite semigroup in which all multiplications, left or right, act as permutations. It is not difficult to see that such a semigroup is in fact a group. Indeed, since all multiplications act as permutations of a finite set, some power of any element acts as an identity; hence the semigroup is a monoid. Since every element permutes the monoid both by right and by left multiplication, every element has both a left and a right inverse, and so the monoid must be a group. All our $S_i^i$ are then finite groups (their identites may not coincide with the identity $1_i$ in the category $I$). Denote by $e_i^i$ the identity of this group structure on $S_i^i$; it is the unique idempotent morphism in $S_i^i$. In fact, $e_i^i$ acts as the identity not only on $S_i^i$, but on all $S_i^j$ by right multiplication and on all $S_j^i$ by left multiplication. This is easily seen from the fact that these actions are permutations and the idempotence of $e_i^i$.   

Now consider the subgraph of the underlying graph of $I$ whose vertices are all the objects of $I$ and whose arrows are those belonging to the sets $S_i^j$. Composition of arrows in $I$ gives this graph a structure of category, with identities $e_i^i$; this follows from the discussion in the previous paragraph. In fact, this category is a groupoid: given $s_i^j\in S_i^j$, take any $s_j^i\in S_j^i$. Then the composition $s_j^is_i^j$ belongs to the group $S_i^i$, so it must be invertible. This means that any morphism $s_i^j$ in our new category is left invertible, so all morphisms are invertible. We will denote this groupoid by $\mathcal G_I$. Notice that it is connected, and the automorphism groups of the vertices are the groups $S_i^i$. In particular, all these groups are isomorphic. We denote this unique finite group by $G_I$. Of course, when regarded as a category with only one object, $G_I$ is equivalent to $\mathcal G_I$. 

The groupoid $\mathcal G_I$ is embedded in $I$ graph-theoretically, but the embedding is not necessarily a functor, since it need not preserve identities. There is, however, a canonical functor $\tau:I\to\mathcal G_I$ which is a left inverse to the embedding of graphs $\mathcal G_I\to I$, and which makes $\mathcal G_I$ into the enveloping groupoid of $I$. We do not require this last fact, but we will define the mentioned functor $\tau$; it is simply the map which acts on morphisms as follows:

$$
\begin{diagram}[labelstyle=\scriptstyle]
i               &\rTo^f                     &j\\
\uTo^{e_i^i}    &\Big\downarrow             &\dTo_{e_j^j}\\
i               &\rTo^{\tau f}              &j  
\end{diagram}
$$
The properties of $e_i^i$ noted above prove that the restriction of $\tau$ to the subgraph $\mathcal G_I\subset I$ is the identity, and also that $\tau$ is indeed a functor. 

The following result will be useful in dealing with the categories ${\bf Set}$ and $_R\mathcal M$ in the next section.

\begin{proposition}\prlabel{cat<>grp}

Let $I$ be a small connected category with an IS consisting of the sets $(S_i^j)$, and let $\mathcal C$ be any complete, cocomplete category. Then $I$ is $\mathcal C$-Frobenius if and only if the group $G_I$ (regarded as a category) is $\mathcal C$-Frobenius.

\end{proposition}

Before embarking on the proof, we need some preparations. Denote by $M$ the two-element monoid $\{1,e\}$, where $1$ is the identity and $e$ is idempotent. Then, regarding $M$ as a one-object category, we have the following simple result:

\begin{lemma}\lelabel{mon2}

$M$ is $\mathcal C$-Frobenius for any complete, cocomplete category $\mathcal C$. 

\end{lemma}

\begin{proof}

A functor $M\to\mathcal C$ is an action of $M$ on some object $c\in\mathcal C$, i.e. a monoid morphism $M\to {\rm Hom}(c,c)$. For such a functor $F$, glue the limting and the colimiting cone into the following commutative diagram:

$$
\begin{diagram}[labelstyle=\scriptstyle]
&                          &                       &c                      &                   &\\
&                          &\ruTo^{\phi}           &\dTo^{Fe}              &\rdTo^{\psi}       &\\
&\varprojlim F             &\rTo^{\phi}            &c                      &\rTo^{\psi}        &\varinjlim F 
\end{diagram}
$$

Because $e$ is idempotent, we get a cone 
$$
\begin{diagram}[labelstyle=\scriptstyle]
&                     &                  &c\\
&                     &\ruTo^{Fe}        &\dTo_{Fe}\\
&c                    &\rTo^{Fe}         &c
\end{diagram}
$$
which induces a unique morphism $\xi:c\to\varprojlim F$ such that $Fe=\phi\xi$. From the uniqueness of $\xi$ we get $\xi\circ Fe=\xi$. Now the commutative diagram
$$
\begin{diagram}[labelstyle=\scriptstyle]
&\varprojlim F             &\rTo^{\phi}            &c               &\rTo^{\xi}                 &\varprojlim\\
&                          &\rdTo_{\phi}           &\dTo^{Fe}       &\ldTo_{\phi}               &\\
&                          &                       &c               &                           & 
\end{diagram}
$$
and the universality of the limit prove that $\xi\phi$ is the identity of $\varprojlim F$. 

Dually, one finds $\eta:\varinjlim F\to c$ through which $Fe$ factors, with the properties $Fe\circ\eta=\eta$ and $\psi\eta=1_{\varinjlim F}$. Putting all of this together we see that the composition $\xi\eta:\varinjlim F\to\varprojlim F$ is the inverse of the natural morphism $\psi\phi:\varprojlim F\to\varinjlim F$ (bottom row of the first diagram above). All the constructions used above are natural with respect to $F$, so we get a natural isomorphism $\varprojlim\cong\varinjlim$, as desired.
\end{proof}

Let $G$ be a semigroup, and denote by $G^+$ the monoid obtained by adjoining an identity to $G$. As a set, it consists of $G$ together with an element $1$; multiplication on $G$ is the one inherited from the semigroup structure of $G$, and $1$ acts as a unit on $G^+=G\cup\{1\}$. When $G$ was a group to begin with (or more generally a monoid), we denote its unit by $e$. Notice that $e$ is an idempotent in $G^+$, but it is no longer the unit for the multiplication in $G^+$. In the proof of \prref{cat<>grp} we make use of the following lemma:

\begin{lemma}\lelabel{G>G+}

Let $\mathcal C$ be any complete and cocomplete category, and let $G$ be a $\mathcal C$-Frobenius group. Then the monoid $G^+$ is also $\mathcal C$-Frobenius.

\end{lemma}

\begin{proof}

The two-element monoid $M$ from the previous lemma is embedded in $G^+$ as $\{1,e\}$, where $1$ is the identity of $G^+$ and $e$ is the identity of $G$. A functor $F:G^+\to\mathcal C$ is an action of the monoid $G^+$ on some object $c\in C$. Restrict this action to the submonoid $M\le G^+$, and let $\phi:d\to c$ be the limiting cone of the restriction $F|_M$. We construct an action $F^*$ of $G$ on $d$ as follows: for every $s\in G$ we have a commutative diagram

$$
\begin{diagram}[labelstyle=\scriptstyle]
&                &                 &c              &\rTo^{Fs}                 &c\\
&                &\ruTo^{\phi}     &\dTo^{Fe}      &                          &\dTo_{Fe}\\
&d               &\rTo^{\phi}      &c              &\rTo^{Fs}                 &c 
\end{diagram}
$$ 
which induces a unique endomorphism $F^*s$ of $d$ making the following diagram commutative:

$$
\begin{diagram}[labelstyle=\scriptstyle]
d            &\rTo^{F^*s}                 &d\\
\dTo^{\phi}  &                            &\dTo_{\phi}\\
c            &\rTo^{Fs}                   &c 
\end{diagram}
$$

That $F^*$ is indeed a functor is easily checked; it must preserve composition by uniqueness because $F$ does, and $F^*e$ is the identity because $\phi:d\to c$ is a cone from $d$ to $F|_M$, and $e$ is a morphism in $M$. 

I claim now that $\varprojlim F$ is naturally isomorphic to $\varprojlim F^*$. Indeed, because $\phi:d\to c$ is limiting, any cone $\varphi:t\stackrel{\cdot}{\rightarrow}c$ which (by definition) makes commutative the diagrams 
$$
\begin{diagram}[labelstyle=\scriptstyle]
&                     &                       &c\\
&                     &\ruTo^{\varphi}        &\dTo_{Fs}\\ 
&t                    &\rTo^{\varphi}         &c 
\end{diagram},\qquad \forall s\in G
$$
must factor through $d$:
$$
\begin{diagram}[labelstyle=\scriptstyle]
&                    &                      &d                 &\rTo^{\phi}             &c\\
&                    &\ruTo                 &\dTo^{F^*s}       &                        &\dTo_{Fs}\\
&t                   &\rTo                  &d                 &\rTo^{\phi}             &c 
\end{diagram},\qquad \forall s\in G
$$

Dually, one constructs an action $F_*$ of $G$ on $\varinjlim F|_M$, and we have a natural isomorphism $\varinjlim F\cong\varinjlim F_*$. Now, because $M$ is always $\mathcal C$-Frobenius (\leref{mon2}), $\varprojlim F|_M\cong\varinjlim F|_M$ naturally. Moreover, recall from the proof of \leref{mon2} that the isomorphism between $\varprojlim F|_M$ and $\varinjlim F|_M$ we have exhibited was precisely the composition of natural maps $\varprojlim F|_M\to c\to\varinjlim F|_M$. The actions $F^*$ and $F_*$ were constructed such that the following diagrams are commutative:
$$
\begin{diagram}[labelstyle=\scriptstyle]
\varprojlim F|_M                   &\rTo                  &c                 &\rTo                &\varinjlim F|_M\\
\dTo^{F^*s}                        &                      &\dTo^{Fs}         &                    &\dTo^{F_*s}\\
\varprojlim F|_M                   &\rTo                  &c                 &\rTo                &\varinjlim F|_M 
\end{diagram}\qquad    \forall s\in G
$$ 
Hence, upon identifying the limit and colimit of $F|_M$ by the given isomorphism, the action $F_*$ is identified to $F^*$. The conclusion now follows from the hypothesis that $G$ is $\mathcal C$-Frobenius.
\end{proof}

Finally, we are ready to prove \prref{cat<>grp}

\renewcommand{\proofname}{Proof of \prref{cat<>grp}}

\begin{proof}

We have noticed in the discussion above that $G_I$ and $\mathcal G_I$ are equivalent categories, so we can replace $G_I$ with $\mathcal G_I$ in the statement of the proposition. 

Assume first that $I$ is $\mathcal C$-Frobenius. Since $\tau:I\to\mathcal G_I$ is a retraction onto the subgraph $\mathcal G_I\to I$, it is bijective on objects and surjective on morphisms. From this it follows immediately that for every $c\in\mathcal C$ the cones $c\stackrel{\cdot}{\rightarrow}F$ coincide with the cones $c\stackrel{\cdot}{\rightarrow}F\tau$. Consequently, the canonical morphism $\varprojlim F\to\varprojlim F\tau$ is an isomorphism. Similarly, $\varinjlim F$ is isomorphic to $\varinjlim F\tau$, naturally in $F$. Applying the $\mathcal C$-Frobenius property to the functors in $\mathcal C^I$ of the form $F\tau$, this discussion implies that $\mathcal G_I$ and hence $G_I$ must be Frobenius as well. 

Conversely, assume that $\mathcal G_I$ (and so $G_I$) is $\mathcal C$-Frobenius. For each object $i\in I$, denote by $M_i$ the submonoid of ${\rm Hom}(i,i)$ consisting of the elements of $S_i^i$ together with the identity. If $S_i^i$ already contains the identity, then $M_i$ is isomorphic to the group $G_I\cong S_i^i$. Otherwise, it will be isomorphic to the monoid denoted above by $G_I^+$. Either way, we know (\leref{G>G+}) that $M_i$ is a $\mathcal C$-Frobenius monoid. 

Given an object $i\in I$ and a functor $F\in\mathcal C^I$, let $F_i$ be the restriction $F|_{M_i}$. If we manage to prove that $\varprojlim F\cong\varprojlim F_i$ naturally (for a fixed $i\in I$), then the dual argument would apply to show that $\varinjlim F\cong\varinjlim F_i$; from the fact that $M_i$ is Frobenius it would then follow that $I$ is also. Hence it remains to prove that there is a natural isomorphism $\varprojlim F\cong\varprojlim F_i$. 

Let $\phi_i:d_i\to F(i)$ be the limiting cones for $F_i$. for objects $i,j\in I$, consider an arbitrary morphism $f_i^j\in {\rm Hom}(i,j)$. I claim that there is a unique morphism $\phi_i^j:d_i\to d_j$ making the following diagram commutative, and that moreover, it does not depend on the morphism $f_i^j$:

$$
\begin{diagram}[labelstyle=\scriptstyle]
d_i                &\rTo^{\phi_i^j}                 &d_j\\
\dTo^{\phi_i}      &                                &\dTo_{\phi_j}\\
F(i)               &\rTo^{Ff_i^j}                   &F(j)   
\end{diagram}
$$
Independence of $f_i^j$ is immediate: since $\phi_i:d_i\to F_i$ is a cone from $d_i$ to $F_i=F|_{M_i}$, we have $Fs_i^i\circ\phi_i=\phi_i$ for every morphism $s_i^i\in S_i^i\subseteq M_i$. Composing to the left with $Ff_i^j$ and using the invariance properties of the IS $(S_i^j)$, we get $Ff_i^j\circ\phi_i=Fs_i^j\circ\phi_i$ for {\it any} $s_i^j\in S_i^j$. The existence of $\phi_i^j$ also follows from this discussion, for it follows that composition to the right with any $Fs_j^j,\ s_j^j\in S_j^j$ fixes $Ff_i^j\circ\phi_i$, so this latter morphism gives a cone $d_i\stackrel{\cdot}{\rightarrow}F_j$. 

From the uniqueness of all $\phi_i^j$ (including the cases $i=j$) it follows that they are isomorphisms; more precisely, for every $i,j\in I,\ \phi_j^i$ is the inverse of $\phi_i^j$. From the universality of $\phi_i:d_i\to F(i)$ and the construction of $\phi_i^j$ it follows that every cone $c\stackrel{\cdot}{\rightarrow}F$ factors through maps $\psi_i:c\to d_i$ making commutative the triangles
$$
\begin{diagram}[labelstyle=\scriptstyle]
c                       &                      &\\
\dTo^{\psi_i}           &\rdTo^{\psi_j}        &\\
d_i                     &\rTo^{\phi_i^j}       &d_j
\end{diagram}\qquad   i,j\in I
$$

Now, since $\phi_i^j$ are isomorphisms, this says that $\varprojlim F$ is naturally isomorphic to $d_i=\varprojlim F_i$ (the constructions appearing above are natural in $F$ once we fix an object $i\in I$). We have thus reached the desired conclusion. 
\end{proof}

\renewcommand{\proofname}{Proof}

\section{Special cases: sets and modules}\selabel{3}

In this section we characterize those small $I$ (not necessarily connected) which are ${\bf Set}$-Frobenius and $_R\mathcal M$-Frobenius for a ring $R$. \seref{1} and \seref{2} will allow us to obtain both necessary and sufficient conditions on $I$ in order that it be Frobenius for these categories. We have already remarked in the discussion on varieties of algebras above that ${\bf Set}$ and $_R\mathcal M$ are admissible categories, so the results in \seref{2} apply in both cases. Remember that all our categories are non-empty.

The following theorem describes the {\bf Set}-Frobenius categories:

\begin{theorem}\thlabel{Set}

A small category $I$ is ${\bf Set}$-Frobenius if and only if it is connected and it has an IS consisting of singletons $S_i^j$. 

\end{theorem}

\begin{proof}

Assume $I$ satisfies the conditions in the statement. Then the group $G_I\cong S_i^i,\ \forall i\in I$ introduced in the discussion before \prref{cat<>grp} is the trivial group. From \prref{cat<>grp} we know that in order to conclude that $I$ is Frobenius, it suffices to check that $G_I$ is. It is clear that the trivial group is $\mathcal C$-Frobenius for any complete, cocomplete category $\mathcal C$, and the proof of this implication is finished. 

Conversely, suppose $I$ is ${\bf Set}$-Frobenius. \leref{discr_conn} (b) then tells us that the set $J$ of connected components of $I$, viewed as a discrete category, must be ${\bf Set}$-Frobenius. The only non-empty ${\bf Set}$-frobenius discrete category is the singleton: notice for instance that the product of a non-empty set and at least one copy of the empty set is empty, whereas the disjoint union of all these sets is non-empty. Hence $J$ is a singleton, i.e. $I$ is connected.

Now \thref{Frob>inv} applies to show that $I$ has an IS consisting of finite non-empty sets $S_i^j$. Now we go once more through the argument in the first paragraph, in reverse: \prref{cat<>grp} says that $I$ is ${\bf Set}$-Frobenius if and only if the finite group $G_I$ is, so we have to prove that the only ${\bf Set}$-Frobenius finite group is the trivial group. 

Functors from $G_I$ to ${\bf Set}$ are actions of $G_I$ on a set. They have easily described limits and colimits: the limiting cone of an action of $G_I$ on the set $c$ is the inclusion of the set of points in $c$ fixed by all elements of $G_I$. The colimiting cone, on the other hand, is the canonical projection of $c$ onto the set of orbits of the action (sending each element onto its orbit). In particular, we see that the colimit of an action on a non-empty set is always a non-empty set, whereas one can always find actions on non-empty sets with no fixed points whenever $G_I$ is non-trivial: simply make $G_I$ act on itself by left multiplication, for example. 
\end{proof}

For $R$-modules, the result reads as follows:

\begin{theorem}\thlabel{R-Mod}

Let $R$ be a ring. A small category $I$ is $_R\mathcal M$-Frobenius if and only if it has finitely many components, each of which has an IS consisting of finite sets $S_i^j$ such that $|S_i^j|$ is invertible in the ring $R$.

\end{theorem}

In the course of the proof we will make use of the following result regarding discrete categories:

\begin{proposition}\prlabel{Iov1}

Let $R$ be a ring, and $J$ a set. The discrete category $J$ is $_R\mathcal M$-Frobenius if and only if it is finite. 

\end{proposition}

\begin{proof}

This is \cite[Theorem 2.7]{Iov1}. In that paper it is both an immediate consequence of the main result \cite[Theorem 1.4]{Iov1}, and proved separately using a finiteness result on Frobenius corings (\cite[Theorem 2.3]{Iov1}; see also \cite[$\S$27]{BW} for definitions and relevant results on Frobenius corings). We give here a different proof, relying on another proposition found in \cite{Iov1}.

On the one hand, it is well-known that finite sets are $_R\mathcal M$-Frobenius. In fact, products and coproducts are canonically isomorphic in any additive category. 

Conversely, assume that $J$ is $_R\mathcal M$-Frobenius. Now \cite[Proposition 1.2]{Iov1} says that the canonical map $\displaystyle\bigoplus_J\to\prod_J$ is a natural isomorphism. Consider the composition
$$
\displaystyle
\begin{diagram}[labelstyle=\scriptstyle]
R           &\rTo             &\prod_JR             &\rTo                &\bigoplus_JR   
\end{diagram}
$$
in which the first arrow is the map with all components equal to the identity on $R$, while the second arrow is the inverse of the canonical isomorphism $\displaystyle\bigoplus_J\to\prod_J$. It is a morphism from $R$ to $\displaystyle\bigoplus_JR$ having the property that the image of $R$ is not contained in any $\displaystyle\bigoplus_{J'}R$ for $J'\subsetneq J$. As $R$ is a finitely generated $R$-module, however, its image is certainly contained in a finite direct sum. Hence $J$ must be finite.
\end{proof}

\begin{remark}\relabel{smobj}

It is clear that the direct sum of infinitely many non-zero modules is strictly smaller than their direct product. However, note that the proof, arranged as above, applies to all (complete, cocomplete) abelian categories having a non-zero {\it small object}. We say that an object $x$ in a category with coproducts is small if any morphism of $x$ to a coproduct factors through a finite coproduct. Indeed, \cite[Proposition 1.2]{Iov1} covers this situation as well (and in fact holds for all categories enriched over the category of commutative monoids and having a zero object), and all we need to do is replace $R$ in the above proof with a small, non-zero object.

A cocomplete abelian category with a small projective generator is equivalent to some $_R\mathcal M$ (\cite[Chapter 4, exercises E and F]{Fr}). There are, however, examples of complete, cocomplete abelian categories with a non-zero small object and which are not equivalent to some $_R\mathcal M$. We give such an example below.
\end{remark}

\begin{example}[\textit{Torsion modules}]\exlabel{1}

Let $R$ be a DVR (discrete valuation ring), and let $\mathcal C$ be the full subcategory of $_R\mathcal M$ consisting of torsion modules. $\mathcal C$ is an abelian category, because kernels, cokernels, finite direct sums, etc. of morphisms of torsion modules are morphisms of torsion modules. Completeness and cocompleteness are, again, easily checked: the direct sum in $_R\mathcal M$ is also the direct sum in $\mathcal C$, and the direct product in $\mathcal C$ is the torsion of the direct product in $_R\mathcal M$. Finally, the category has non-zero small objects: any non-zero finitely generated torsion module will do. A small projective object in $\mathcal C$ must be finitely generated, and the structure theorem for finitely generated modules over a PID now easily shows that $\mathcal C$ has no non-zero small projectives, hence cannot be equivalent to some $_S\mathcal M$. 
\end{example}

\begin{remark}\relabel{genex}

Although we will not prove this here, with a little more work, it can be shown that the previous example still works if $R$ is taken to be any noetherian local integral domain (which is not a field). 
\end{remark}

At the other end of the spectrum, when working with connected categories, we will need the following characterization of $_R\mathcal M$-Frobenius groups:

\begin{proposition}\prlabel{Frobgr}

Let $R$ be a ring and $G$ a group, regarded as a one-object category. $G$ is $_R\mathcal M$-Frobenius if and only if it is finite, and the natural number $|G|$ is invertible in $R$. 

\end{proposition}

\begin{proof}

Functors $G\to\ _R\mathcal M$ are precisely $R$-modules with a $G$ action, or, in other words, $R[G]$-modules. The diagonal functor $_R\mathcal M\to(_R\mathcal M)^G$ associates to each $R$-module the same module with trivial $G$ action. This means that one can identify the diagonal functor with the restriction of scalars from $R$ to $R[G]$ through the augmentation $\varepsilon:R[G]\to R$ (the unique ring morphism sending each element of $G\subset R[G]$ to the identity $1_R\in R$). 

The problem has now been reduced to the classical question of deciding when a restriction of scalars is Frobenius. By a well-known result of Morita (\cite{Mo} or \cite[Theorem 3.15]{Men}), restriction of scalars through a ring morphism $A\to B$ is Frobenius if and only if $B$ is left $A$-projective and finitely generated, and $B\cong\ _A{\rm Hom}(B,A)$ as $(B,A)$-bimodules. We are going to apply this characterization to the ring extension $\varepsilon:R[G]\to R$. 

$R$ is left $R[G]$-projective if and only if the augmentation $\varepsilon:R[G]\to R$ splits through some left $R[G]$-module map $\eta:R\to R[G]$. For any such splitting, $\eta(1)$ is some element $\displaystyle\sum_{g\in G}a_gg$ of $R[G]$ fixed by left multiplication with any element of $G$. This shows at once that $G$ must be finite, and that $\displaystyle\eta(1)=a\sum_{g\in G}g$. Finally, from $\varepsilon\circ\eta={\rm id}_R$ we find that $a\in R$ must in fact be the inverse of $|G|$. Conversely, if $|G|<\infty$ is invertible in $R$, simply consider the $R$-module map sending $1\in R$ to $\displaystyle |G|^{-1}\sum_{g\in G}g\in R[G]$; clearly, it is a splitting for $\varepsilon$. 

We still have to prove that when $G$ satisfies the condtitions in the statement of the proposition (and hence, as we have just seen, $R$ is left $R[G]$-projective), we also have an isomorphism $\displaystyle R\cong\ _{R[G]}{\rm Hom}(R,R[G])$ in $_R\mathcal M_{R[G]}$. The second term is canonically isomorphic to the $(R,R[G])$-sub-bimodule of $R[G]$ generated by the central idempotent $e=\displaystyle|G|^{-1}\sum_Gg$; there is an obvious $(R,R[G])$-bimodule isomorphism of $R$ onto this bimodule, sending $1$ to $e$.  
\end{proof}

We are now ready to prove the theorem.

\renewcommand{\proofname}{Proof of \thref{R-Mod}}
\begin{proof}

Since $_R\mathcal M$ is a complete, cocomplete category with a zero object, points (a) and (c) of \leref{discr_conn} show that $I$ is Frobenius if and only if (i) its set of connected components $J$ is Frobenius, and (ii) each connected component is Frobenius. Hence the problem breaks up into the discrete and the connected case. 

\prref{Iov1} says that the component set is $_R\mathcal M$-Frobenius if and only if it is finite. In the connected case we can apply the results in \seref{2}. \thref{Frob>inv} and \prref{cat<>grp} together imply that a connected category is $_R\mathcal M$-Frobenius if and only if it has an IS such that the group $G_I$ is $_R\mathcal M$-Frobenius. Finally, apply \prref{Frobgr} to finish the proof.
\end{proof}
\renewcommand{\proofname}{Proof}

\section{Some open problems}\selabel{4}

The problem posed here, of finding the $\mathcal C$-Frobenius categories $I$ for a fixed complete and cocomplete $\mathcal C$, has variations which would make interesting topics for further inquiry. We give only a few examples. 

For one thing, we would like to extend the results obtained in this paper to various categories (or perhaps large classes of categories) which were not covered here. One conspicuous example is that of the category of (left or right) comodules over some $R$-coring $C$. This would cover the case of $R$-modules, since these are the simply the comodules over the Sweedler coring $R$ over $R$ (\cite[Examples 17.3 and 18.5]{BW}). Choose right comodules, in order to fix the notation. Because we want the category $\mathcal M^C$ of right comodules to be complete and cocomplete, we impose the condition that $_RC$ be flat (see \cite[Theorem 18.13]{BW}).

\begin{problem}

Given a ring $R$ and an $R$-coring $C$ which is flat as a left $R$-module, find the $\mathcal M^C$-Frobenius small categories $I$. 

\end{problem}

Even within the realm of admissible categories, treated here, the results we have proven give rise to some interesting questions. For example, \thref{Frob>inv} and \prref{cat<>grp} together reduce the problem of finding the connected $\mathcal C$-Frobenius categories to that of finding the $\mathcal C$-Frobenius finite groups, whenever $\mathcal C$ is admissible. We have already seen two classes of groups arising as the class of $\mathcal C$-Frobenius finite groups for various $\mathcal C$: the trivial group if $\mathcal C={\bf Set}$, and the finite groups whose cardinality is invertible in $R$ for $\mathcal C=\ _R\mathcal M$. Can all such classes of finite groups be described?

\begin{problem}

Which classes of finite groups arise as the class of $\mathcal C$-Frobenius finite groups for some admissible category $\mathcal C$? 

\end{problem}

We can turn this question around, and ask for a characterization of those admissible categories $\mathcal C$ having the property that the only $\mathcal C$-Frobenius finite group is the trivial group. We have already seen in \thref{Set} that ${\bf Set}$ is such a category. Although we do not prove this here, it is not difficult to see that ${\bf Grp}$, the category of groups, is another example. Note that {\bf Grp} is a variety of algebras, so it is indeed admissible.

\begin{problem}

Find simple necessary and sufficient (or, alternatively, only sufficient) conditions on an admissible category $\mathcal C$ in order that the only $\mathcal C$-Frobenius finite group be the trivial group. 

\end{problem} 

\section*{Aknowledgement}

The author wishes to thank Professor Gigel Militaru, who posed the problem and suggested this line of inquiry, for the insight gained through countless discussions on the topic, as well as the referee for valuable suggestions on how to revise an initial version of this paper. 



\end{document}